\documentclass[final,leqno,onefignum,onetabnum]{siamltex1213}

\usepackage{amsfonts,amsmath,amssymb}
\usepackage{graphicx}
\usepackage[dvips]{color}

\usepackage{algorithm}
\usepackage{algorithmic}

\begin{document}

\newcommand{\NN}{\mathbb{N}}
\newcommand{\ZZ}{\mathbb{Z}}
\newcommand{\RR}{\mathbb{R}}
\newcommand{\CC}{\mathbb{C}}

\title{A trivariate interpolation algorithm using a cube-partition searching procedure
} 

\author{Roberto Cavoretto\footnotemark[2]\ \footnotemark[3]
\and Alessandra De Rossi\footnotemark[2]\ \footnotemark[4]
}

\renewcommand{\thefootnote}{\fnsymbol{footnote}}

\footnotetext[2]{Department of Mathematics \lq\lq G. Peano\rq\rq, University of Torino, via Carlo Alberto 10, I--10123 Torino, Italy ({\tt roberto.cavoretto@unito.it, alessandra.derossi@unito.it}).}
\footnotetext[3]{This author's work was partially supported by the University of Torino via grant \lq\lq Approssimazione di dati sparsi e sue applicazioni\rq\rq.}
\footnotetext[4]{This author's work was partially supported by the GNCS-INdAM.}

\maketitle
\slugger{mms}{xxxx}{xx}{x}{x--x}

\begin{abstract}
In this paper we propose a fast algorithm for trivariate interpolation, which is based on the partition of unity method for constructing a global interpolant by blending local radial basis function interpolants and using locally supported weight functions. The partition of unity algorithm is efficiently implemented and optimized by connecting the method with an effective cube-partition searching procedure. More precisely, we construct a cube structure, which partitions the domain and strictly depends on the size of its subdomains, so that the new searching procedure and, accordingly, the resulting algorithm enable us to efficiently deal with a large number of nodes. Complexity analysis and numerical experiments show high efficiency and accuracy of the proposed interpolation algorithm.
\end{abstract}



\pagestyle{myheadings}
\thispagestyle{plain}
\markboth{ROBERTO CAVORETTO AND ALESSANDRA DE ROSSI}{3D INTERPOLATION USING A CUBE-PARTITION SEARCHING PROCEDURE}

\section{Introduction}
The problem of constructing fast algorithms for multivariate approximation of scattered data points has recently interested many researchers, who work in various areas of applied mathematics and scientific computing such as interpolation, approximation theory, neural networks, computer aided geometric design (CAGD) and machine learning, to name a few. So we often need to have numerical algorithms, which allow us to efficiently deal with a large number of points, not only in one or two dimensions but also in higher dimensions, as it usually occurs in several applications (see, e.g., \cite{Fasshauer07,Wendland05} and references therein). 

Though there exist several numerical algorithms and alternative techniques for bivariate interpolation to scattered data, the problem of efficiently approximating many thousands or millions of three dimensional data does not seem to be much considered in the literature, with the exception of a few cases such as in \cite{Beatson00,Cherrie02,Lazzaro02,Renka88c,Thacker10}; a comparison of radial basis function (RBF) methods in the 3D setting can be found in \cite{Bozzini02}.

Since mesh-based methods require some sort of an underlying computational mesh, i.e. any triangulation of the domain, their construction is a rather difficult task, already in two dimensions, where the mesh generation turns out usually to be one of the most time consuming part. For this reason, in the following we focus on a \textsl{meshfree} or \textsl{meshless} approximation. More precisely, here we consider the partition of unity method, which involves the use of RBFs as local approximants and of locally supported weight functions (see \cite{Wendland02}). Further details on the origin of the partition of unity method can be found in \cite{Babuska97,Melenk96}. Moreover, some other examples of local approaches involving modified Shepard's methods and different searching procedures can be found in \cite{Allasia11, Berry99, Cavoretto14c, Iyer06, Lazzaro02, Renka88a, Renka88c, Thacker10}.

Starting from the previous work \cite{Cavoretto14a}, where an efficient algorithm with a new cell-based searching procedure is presented for bivariate interpolation of large scattered data sets, in this paper we directly extend it to trivariate case, obtaining in this way a new fast algorithm for interpolation, which can briefly be summarized in three stages as follows:
\begin{enumerate}
	\item[(i)] partition the domain into a suitable number of cubes; 
	\item[(ii)] consider an optimized cube-partition searching procedure establishing the minimal number of cubes to be examined, in order to localize the subset of nodes belonging to each subdomain; 
	\item[(iii)] apply the partition of unity method combined with local RBFs. 
\end{enumerate}
In particular, the algorithm is characterized by the construction of a \textsl{cube-partition searching procedure}, whose origin comes from the repeated use of a \textsl{quicksort} routine with respect to different directions, which enables us to pass from unordered to ordered data structures. Moreover, this technique is strictly related to the construction of a partition of the domain in cubes and depends on the size of its subdomains, thus producing a nearest neighbor searching procedure, which is particularly efficient in local interpolation methods. Numerical experiments show efficiency and accuracy of the cube algorithm.

The paper is organized as follows. In Section \ref{PUM} we recall some theoretical results, giving a general description of the partition of unity method, which makes use of RBFs as local approximants. In Section \ref{PUM_ALG}, we present in detail the cube-partition algorithm for trivariate interpolation, which is efficiently implemented and optimized by using a nearest neighbor searching procedure. Computational complexity and storage requirements of the interpolation algorithm are analyzed as well. In Section \ref{num_res}, we show numerical results concerning efficiency and accuracy of the partition of unity algorithm. Finally, Section \ref{concl} deals with conclusions and future work.


\section{Partition of unity interpolation} \label{PUM}
Let ${\cal X}_n=\{\boldsymbol{x}_i, i = 1,2, \ldots , n \}$ be a set of distinct data points or nodes, arbitrarily distributed in a domain $\Omega \subseteq \RR^N$, $N \geq 1$, with an associated set ${\cal F}_n=\{f_i, i = 1,2,  \ldots , n \}$ of data values or function values, which are obtained by sampling some (unknown) function $f:\Omega \rightarrow \RR$ at the nodes, i.e., $f_i=f(\boldsymbol{x}_i)$, $i=1,2,\ldots,n$.

The basic idea of the partition of unity interpolation is to start with a partition of the open and bounded domain $\Omega \subseteq \RR^N$ into $d$  subdomains $\Omega_j$ such that $\Omega \subseteq \bigcup_{j=1}^{d} \Omega_j$ with some mild overlap among the subdomains. Associated with these subdomains we choose a partition of unity, i.e. a family of compactly supported, non-negative, continuous functions $W_j$ with $\text{supp}(W_j) \subseteq \Omega_j$ such that 
\begin{eqnarray} \label{pu_f}
\sum_{j=1}^{d} W_j(\boldsymbol{x}) = 1.
\end{eqnarray} 
For each subdomain $\Omega_j$ we consider a local approximant $R_j$ and form then the global approximant
\begin{eqnarray}
\label{pui}
	{\cal I}(\boldsymbol{x})= \sum_{j=1}^{d} R_j(\boldsymbol{x}) W_j(\boldsymbol{x}), \hspace{1cm} \boldsymbol{x} \in \Omega.
\end{eqnarray}
Here $R_j:\Omega \rightarrow \mathbb{R}$ defines a RBF interpolant of the form
\begin{align*}
R_j(\boldsymbol{x})= \sum_{k=1}^{\bar{n}_j} c_k \phi (\left\| \boldsymbol{x} - \boldsymbol{x}_k \right\|_2),
\label{rad1}
\end{align*}
where $ \phi: [0, \infty) \rightarrow \mathbb{R}$ represents a radial basis function, $||  \cdot  ||_2$ denotes the Euclidean norm, and $\bar{n}_j$ indicates the number of data points in $\Omega_j$. Furthermore, $R_j$ satisfies the interpolation conditions
\begin{equation} \label{int1}
R_j( \boldsymbol{x}_i)=f_i, \quad i=1,2, \ldots, \bar{n}_j.
\end{equation}
Note that if the local approximants satisfy the interpolation conditions \eqref{int1}, then the global approximant also interpolates at this node, i.e. 
\begin{equation}
{\cal I}(\boldsymbol{x}_i)=f_i, \quad i=1,2,\ldots,\bar{n}_j. \nonumber
\end{equation}

Solving the $j$-th interpolation problem \eqref{int1} leads to a system of linear equations of the form
    \begin{align*}
\left[
\begin{array}{cccc}
\phi (||\boldsymbol{x}_1-\boldsymbol{x}_1||_2)   & \phi (||\boldsymbol{x}_1-\boldsymbol{x}_2||_2) & \cdots  & \phi (||\boldsymbol{x}_1-\boldsymbol{x}_{\bar{n}_j}||_2)        \\
\phi (||\boldsymbol{x}_2-\boldsymbol{x}_1||_2)   & \phi (||\boldsymbol{x}_2-\boldsymbol{x}_2||_2) & \cdots  & \phi (||\boldsymbol{x}_2-\boldsymbol{x}_{\bar{n}_j}||_2)        \\
\vdots         & \vdots        & \vdots      & \vdots        \\
\phi (||\boldsymbol{x}_{\bar{n}_j}-\boldsymbol{x}_1||_2)   & \phi (||\boldsymbol{x}_{\bar{n}_j}-\boldsymbol{x}_2||_2) & \cdots  & \phi (||\boldsymbol{x}_{\bar{n}_j}-\boldsymbol{x}_{\bar{n}_j}||_2)        
\end{array}
\right]
\left[
\begin{array}{c}
c_1\\
c_2\\
\vdots\\
c_{\bar{n}_j}
\end{array}
\right]
=
\left[
\begin{array}{c}
f_1\\
f_2\\
\vdots\\
f_{\bar{n}_j}
\end{array}
\right],
  \end{align*}
or simply
\begin{equation} \label{mat}
 \Phi\boldsymbol{c}=\boldsymbol{f}. \nonumber
\end{equation}
In particular, the interpolation problem is well-posed, i.e., a solution to the problem exists and is unique, if and only if the matrix $\Phi$ is nonsingular. A sufficient condition to have nonsingularity is that the corresponding matrix is positive definite. In fact, if the matrix $\Phi$ is positive definite, then all its eigenvalues are positive and therefore $\Phi$ is nonsingular (see, e.g., \cite{Fasshauer07}).

Though the theory of RBFs is here considered, for brevity we do not report basic definitions and theorems, referring to \cite{Buhmann03, Fasshauer07, Iske11, Wendland05} for a more detailed analysis. Then, we give the following definition (see \cite{Wendland02}).

\begin{definition}
Let $\Omega \subseteq \RR^N$ be a bounded set. Let $\{\Omega_j\}_{j=1}^{d}$ be an open and bounded covering of $\Omega$. This means that all $\Omega_j$ are open and bounded and that $\Omega \subseteq \bigcup_{j=1}^{d} \Omega_j$. Set $\delta_j = {\rm diam}(\Omega_j)=\sup_{\boldsymbol{x},\boldsymbol{y} \in \Omega_{j}} ||\boldsymbol{x}-\boldsymbol{y}||_2$. We call a family of nonnegative functions $\{W_j\}_{j=1}^{d}$ with $W_j \in C^k(\RR^N)$ a $k$-stable partition of unity with respect to the covering $\{\Omega_j\}_{j=1}^{d}$ if
\begin{enumerate}
	\item[1)] ${\rm supp}(W_j) \subseteq \Omega_j$;
	\item[2)] $\sum_{j=1}^{d} W_j(\boldsymbol{x}) \equiv 1$ on $\Omega$;
	\item[3)] for every $\beta \in \NN_0^N$ with $|\beta| \leq k$ there exists a constant $C_{\beta} > 0$ such that
\begin{equation}
	||D^{\beta}W_j||_{L_{\infty}(\Omega_j)}\leq C_{\beta}/\delta_j^{|\beta|}, \nonumber
\end{equation}
for all $1\leq j \leq d$.
\end{enumerate}
\end{definition}

In agreement with the statements in \cite{Wendland02}, we require additional regularity assumptions on the \textsl{covering} $\{\Omega_j\}_{j=1}^{d}$.

\begin{definition} \label{defpr}
Suppose that $\Omega \subseteq  \RR^N$ is bounded and ${\cal X}_n=\left\{\boldsymbol{x}_i, i=1,2,\ldots,n\right\} \subseteq \Omega$ are given. An open and bounded covering $\{\Omega_j\}_{j=1}^{d}$ is called regular for $(\Omega,{\cal X}_n)$ if the following properties are satisfied:
\begin{itemize}
	\item[(a)] for each $\boldsymbol{x} \in \Omega$, the number of subdomains $\Omega_j$ with $\boldsymbol{x} \in \Omega_j$ is bounded by a global constant $K$;
	\item[(b)] each subdomain $\Omega_j$ satisfies an interior cone condition;
	\item[(c)] the local fill distances $h_{{\cal X}_{j}, \Omega_j}$, where ${\cal X}_{j}={\cal X}_n \cap \Omega_j$, are uniformly bounded by the global fill distance $h_{{\cal X}_n, \Omega}$, i.e.
	\begin{align*} 
		h_{{\cal X}_{n}, \Omega} = \sup_{\boldsymbol{x} \in \Omega}\min_{\boldsymbol{x}_k\in {\cal X}_n} ||\boldsymbol{x}-\boldsymbol{x}_k||_2. 
	\end{align*}
\end{itemize}
\end{definition}

Property (a) is required to ensure that the sum in \eqref{pui} is actually a sum over at most $K$ summands. Since $K$ is independent of $n$, unlike $d$, which should be proportional to $n$, this is essential to avoid losing convergence orders. It is crucial for an efficient evaluation of the global interpolant that only a constant number of local approximants has to be evaluated. In such way, it should be possible to locate those $K$ indices in constant time. Properties (b) and (c) are important for employing the estimates on RBF interpolants (see \cite{Wendland05}).

Moreover, we are able to formulate the following theorem, which yields the polynomial precision and controls the growth of error estimates, denoting by $\pi_s^N := \pi_s(\RR^N)$ the set of polynomials of degree at most $s$ (see, e.g., \cite{Wendland05}).

\begin{theorem} \label{TP}
	Suppose that $\Omega \subseteq \RR^N$ is compact and satisfies an interior cone condition with angle $\theta \in (0,\pi/2)$ and radius $r > 0$. Let $s \in \NN$ be fixed and there exist constants $h_0, C_1, C_2 > 0$ depending only on $N, \theta, r$ such that $h_{{\cal X}_n, \Omega} \leq h_{0}$. Then, for all ${\cal X}_n=\{\boldsymbol{x}_i, i=1,2,\ldots,n\} \subseteq \Omega$ and all $\boldsymbol{x} \in \Omega$, there exist functions $u_k:\Omega \rightarrow \RR$, $k=1,2,\ldots,n$, such that
\begin{enumerate}
	\item[(1)] $\sum_{k=1}^{n} u_k(\boldsymbol{x})p(\boldsymbol{x}_k) = p(\boldsymbol{x})$, \ \mbox{ for all } $p \in \pi_{s}(\RR^N)$; 
	\item[(2)] $\sum_{k=1}^{n} \left|u_k(\boldsymbol{x})\right|\leq C_1$;
	\item[(3)] $u_j(\boldsymbol{x}) = 0$ provided that $||\boldsymbol{x}-\boldsymbol{x}_j||_2 > C_2 h_{{\cal X}_n, \Omega}$.
\end{enumerate}
\end{theorem}

Therefore, after defining the space $C_{\nu}^k(\RR^N)$ of all functions $f \in C^k$ whose derivatives of order $|\beta|=k$ satisfy $D^{\beta}f(\boldsymbol{x})= {\cal O}(||\boldsymbol{x}||_2^{\nu})$ for $||\boldsymbol{x}||_2 \rightarrow 0$, we consider the following convergence result (see, e.g., \cite{Fasshauer07, Wendland05}).

\begin{theorem}
	Let $\Omega \subseteq  \RR^N$ be open and bounded and suppose that ${\cal X}_n = \{\boldsymbol{x}_i, i=1,$ $2,\ldots,n \}\subseteq \Omega$. Let $\phi \in C_{\nu}^k(\RR^N)$ be a strictly conditionally positive definite function of order $m$. Let $\{\Omega_j\}_{j=1}^{d}$ be a regular covering for $(\Omega, {\cal X}_n)$ and let $\{W_j\}_{j=1}^{d}$ be $k$-stable for $\{\Omega_j\}_{j=1}^{d}$. Then the error between $f \in {\cal N}_{\phi}(\Omega)$, where ${\cal N}_{\phi}$ is the native space of $\phi$, and its partition of unity interpolant (\ref{pui}) can be bounded by
\begin{equation}
	|D^{\beta}f(\boldsymbol{x}) - D^{\beta}{\cal I}(\boldsymbol{x})| \leq C h_{{\cal X}_n, \Omega}^{(k+\nu)/2 - |\beta|} |f|_{{\cal N}_{\phi}(\Omega)}, \nonumber
\end{equation}
for all $\boldsymbol{x} \in \Omega$ and all $|\beta| \leq k/2$. 
\end{theorem}  

Comparing this convergence result with the global error estimates (see e.g. \cite{Wendland05}), we note that the partition of unity preserves the local approximation order for the global fit. This means that we can efficiently compute large RBF interpolants by solving small RBF interpolation problems (in parallel as well) and then glue them together with the global partition of unity $\{W_j\}_{j=1}^{d}$. In other words, the partition of unity approach is a simple and effective technique to decompose a large problem into many small problems while at the same time ensuring that the accuracy obtained for the local fits is carried over to the global one. In particular, the partition of unity method can be thought as a Shepard's type interpolation with higher-order data, since local approximations $R_j$ instead of data values $f_j$ are used. 

Finally, we remark that, among several weight functions $\bar{W}_j(\mathbf{x})$ in \eqref{pui}, a possible choice is given by Shepard's weight
\begin{eqnarray} \label{sh_w}
	W_j(\boldsymbol{x}) = \frac{\bar{W}_j(\boldsymbol{x})}{\sum_{k=1}^{d} \bar{W}_k(\boldsymbol{x})}, \hspace{0.5cm} j=1,2,\ldots,d,
\end{eqnarray}
where $\bar{W}_j$ is the inverse of the Euclidean norm $\Vert \cdot \Vert_2$. It constitutes a partition of unity as in \eqref{pu_f}.

\section{Cube-partition algorithm} \label{PUM_ALG}

In this section we propose a new algorithm for trivariate interpolation of large scattered data sets lying on the domain $\Omega = [0,1]^3 \subset \RR^3$. This algorithm, which is based on the partition of unity method for constructing a global interpolant by blending RBFs as local approximants and using locally supported weight functions, is efficiently implemented and optimized by connecting the interpolation method with an effective cube-partition searching procedure. More precisely, the considered approach is characterized by the construction of a \textsl{cube-based structure}, which partitions the domain $\Omega$ in cubes and strictly depends on the dimension of its subdomains. This technique is a direct extension in three-dimensional case of the square-partition searching procedure presented in \cite{Cavoretto14a} for bivariate interpolation, which we briefly recall in Subsection \ref{review}. 

Note that the paper \cite{Cavoretto14a} follows preceding works, where efficient searching procedures based on the partition of the domain in strips or spherical zones are considered (see \cite{Allasia11,Cavoretto10a,Cavoretto12a,Cavoretto14a}). 


\subsection{Review of the 2D square-partition searching procedure} \label{review}

\begin{figure}[ht!]
\begin{center}
\includegraphics[width=16.cm]{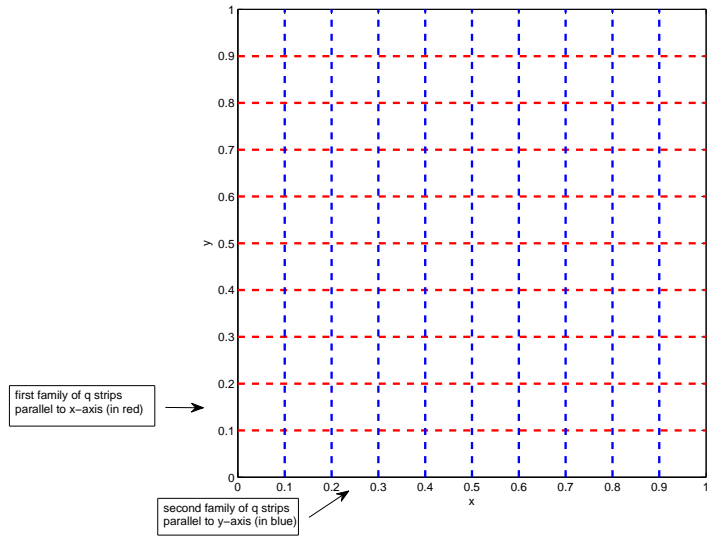}
\end{center}
\caption{Example of orthogonal families of strips.}
\label{double-strips}
\end{figure}

\begin{figure}[ht!]
\begin{center}
\includegraphics[width=16.cm]{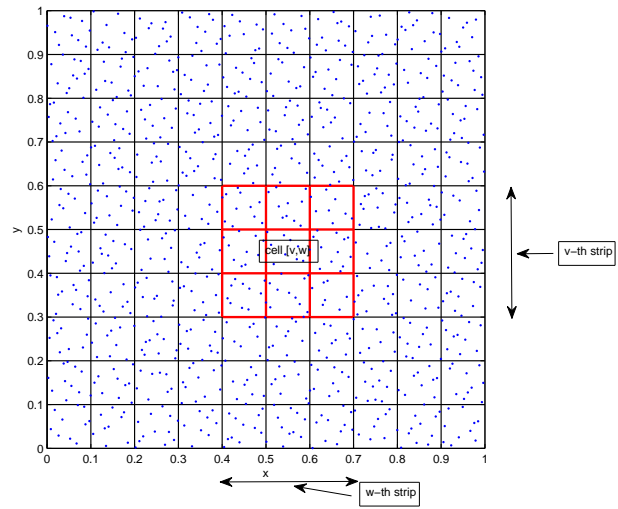}
\end{center}
\caption{Example of square-based structure with a set of scattered data points.}
\label{cell-struct}
\end{figure}

The construction of the 2D searching procedure described in \cite{Cavoretto14a} is obtained by making a partition of the bivariate domain in square cells. They are achieved generating two orthogonal families of parallel strips (see Figure \ref{double-strips}). This approach is combinated with the repeated use of a \textsl{quicksort} routine with respect to different directions. At first, we make a sorting along the $y$-axis on all the points, constructing then a first family of strips parallel to the $x$-axis. Afterwards, we order the points contained in each strip with respect to the $x$-axis direction, and finally we build the second family of strips parallel to the $y$-axis. The outcome is a square-based structure, which allows us to pass from unordered to ordered data structures. Following this idea, we can suitably split up the original data set in ordered and well-organized data subsets. More precisely, we may act as follows: 
\begin{enumerate}
	\item[i)] organize all the data by means of a \textsl{quicksort$_y$ procedure} applied along the $y$-axis (the subscript denotes the sorting direction);
	\item[ii)] consider a first family of $q$ strips, parallel to the $x$-axis and order the points of each strip by using a \textsl{quicksort$_x$ procedure};
	\item[iii)] create a second family of $q$ strips, parallel to the $y$-axis, which orthogonally intersect the first strip family, thus producing a partition of the bivariate domain in square cells (see Figure \ref{cell-struct}).
\end{enumerate}
Note that a specific square cell $k$ is denoted by a double index notation in square brackets, i.e. $k=[v,w]$. 

In order to obtain an efficient searching technique in the localization of points, we connect the interpolation me\-thod with the square-based partition structure, exploiting the data structure and the domain partition previously considered. This result is obtained assuming that the square side is equal to the subdomain radius. Though this choice might seem to be trivial, in practice such an imposition means that the search of the nearby points, an essential aspect of local methods as the partition of unity method, is limited at most to nine squares: the square on which the considered point lies, and the eight neighbouring squares (see Figures \ref{double-strips}--\ref{cell-struct}). The combination between square cell and subdomain sizes constitutes an \textsl{optimal} choice, since it allows us to search the closest points only considering a very small number of them, i.e. taking those points belonging to one of the nine square cells and \textsl{a priori} ignoring all the other ones. Finally, for all those points belonging to the first and last square cells, namely the ones located on or close to the boundary of the domain, we reduce the total number of square cells to be examined.

\subsection{The 3D cube-partition searching procedure} \label{cube-proc}

\begin{figure}[ht!]
\begin{center}
\includegraphics[width=12.cm]{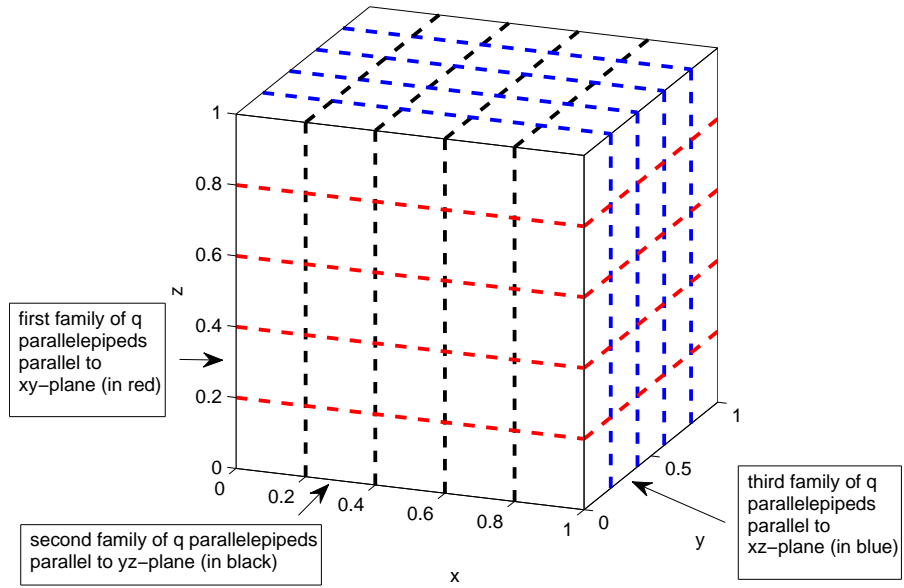}
\end{center}
\caption{Example of orthogonal families of parallelepipeds.}
\label{cube_1}
\end{figure}

As in the 2D case, the basic idea in constructing the 3D searching procedure comes from the repeated use of a \textsl{quicksort} routine with respect to (three) different directions, i.e. along the $z$-axis, the $y$-axis and the $x$-axis, enabling us to pass from unordered to ordered data structures. This process is strictly related to the construction of a partition of the domain, here the unit cube, in smaller cubes. They are obtained generating three orthogonal families of parallelepipeds, while at the same time the original data set is suitably split up in ordered and well-organized data subsets. More precisely, in order to obtain the cube-based structure and then the resulting searching procedure, we may act as follows: 
\begin{enumerate}
	\item[i)] organize all the data by means of a \textsl{quicksort$_z$ procedure} applied along the $z$-axis;
	\item[ii)] consider a first family of $q$ parallelepipeds, parallel to the $xy$-plane, and order the points of each parallelepiped by using a \textsl{quicksort$_x$ procedure}; 
	\item[iii)] create a second family of $q$ parallelepipeds, parallel to the $yz$-plane, which orthogonally intesect the first family, and order the points of each parallelepiped by using a \textsl{quicksort$_y$ procedure};
	\item[iv)] construct a third family of $q$ parallelepipeds, parallel to the $xz$-plane, which orthogonally intesect the two previous families, thus producing a partition of $\Omega$ in cubes (see Figure \ref{cube_1}). 
\end{enumerate}

\begin{figure}[ht!]
\begin{center}
\includegraphics[width=12.cm]{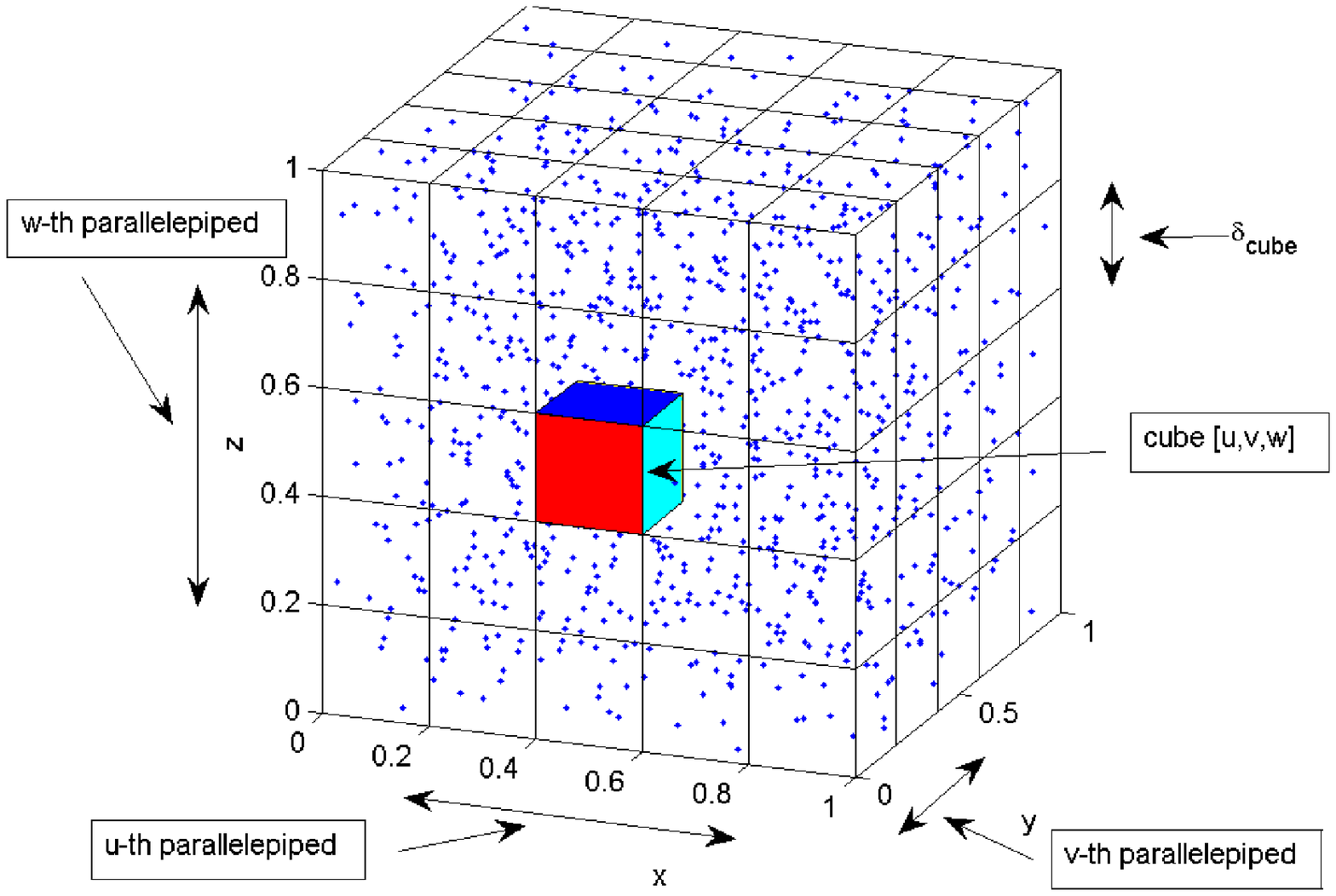}
\end{center}
\caption{Example of cube-based structure with a set of scattered data points.}
\label{cube_2}
\end{figure}

Now, exploiting the data structure and the domain partition, we construct an efficient searching technique to be used in the localization of points, effectively connecting the partition of unity scheme with the cube-partition structure. This result is got assuming that the cube side $\delta_{cube}$ is equal to the subdomain radius $\delta_{subdom}$, i.e. taking $\delta_{cube} \equiv \delta_{subdom}$. From this assumption it follows that the search of the nearby points is limited at most to twenty-seven ($3^3$) cubes: the cube on which the considered point lies, and the twenty-six neighboring cubes (see Figure \ref{cube_2}). From now on, to locate a specific cube $k$, we define a triple index notation using square brackets, i.e. $k=[u,v,w]$, $u,v,w=1,2,\ldots,q$. 

We note that the combination between cube and subdomain sizes provides an \textsl{optimal} choice, since it allows us to search the closest points only considering a very small number of them (that is only those points belonging to one of the twenty-seven cubes) and \textsl{a priori} ignoring all the other points of $\Omega$. Obviously, then, for all those points belonging to cubes close to the boundary of $\Omega$, it will be required a reduction of the total number of cubes to be examined. Further details on this searching procedure are contained in Subsection \ref{cube-algor}, where we give a detailed description of the proposed algorithm.


\subsection{Cube algorithm} \label{cube-algor}

\noindent INPUT: $n$, number of data; ${\cal X}_n=\{(x_i,y_i,z_i), i=1,2,\ldots,n\}$, set of data points; ${\cal F}_n=\{f_i, i=1,2,\ldots,n\}$, set of data values; $d$, number of subdomains; ${\cal C}_d=\{(\bar{x}_i, \bar{y}_i, \bar{z}_i), i=1,2,\ldots,d\}$, set of subdomain points (centres); $s$, number of evaluation points; ${\cal E}_s=\{(\tilde{x}_i,\tilde{y}_i, \tilde{z}_i), i=1,2,\ldots,s\}$, set of evaluation points.
\vskip 3pt
\noindent OUTPUT: ${\cal A}_s=\{{\cal I}(\tilde{x}_i,\tilde{y}_i,\tilde{z}_i), i=1,2,\ldots,s\}$, set of approximated values.
\vskip 3pt
\noindent {\tt Stage 1.} The set ${\cal X}_n$ of nodes and the set ${\cal E}_s$ of evaluation points are ordered with respect to a common direction (e.g. the $z$-axis), by applying a \textsl{quicksort$_z$ procedure}.
\vskip 3pt
\noindent {\tt Stage 2.} For each subdomain point $(\bar{x}_i,\bar{y}_i,\bar{z}_i)$, $i=1,2,\ldots,d$, a local spherical subdomain is constructed, whose spherical radius depends on the subdomain number $d$, i.e.
\begin{eqnarray}
\label{delta}
 \delta_{subdom} = \frac{\sqrt{2}}{\sqrt[3]{d}}.
 \end{eqnarray}
Although other choices $\delta_{subdom}$ are possible, this value is suitably chosen, supposing to have a nearly uniform node distribution and assuming that the ratio $n/d \approx 2^3$.
\vskip 3pt 
\noindent {\tt Stage 3.} A triple structure of intersecting parallelepipeds is constructed as follows:
\begin{enumerate}
	\item[i)] a first family of $q$ parallelepipeds, parallel to the $xy$-plane, is considered taking
\begin{eqnarray} \label{q_par}
	q = \left\lceil \frac{1}{\delta_{subdom}} \right\rceil,
\end{eqnarray}
and a \textsl{quicksort$_x$ procedure} is applied to order the nodes belonging to each parallelepiped;
	\item[ii)] a second family of $q$ parallelepipeds, parallel to the $yz$-plane, is constructed and a \textsl{quicksort$_y$ procedure} is used to order the nodes belonging to each of the resulting parallelepipeds;
	\item[iii)] a third family of $q$ parallelepipeds, parallel to the $xz$-plane, is considered.
\end{enumerate} 
Note that each of the three families of parallelepipeds are ordered and numbered from 1 to $q$; the choice in (\ref{q_par}) follows directly from the side length of the domain, i.e. the unit cube, and the subdomain radius $\delta_{subdom}$.
\vskip 3pt 
\noindent {\tt Stage 4.} The unit cube is partitioned by a cube-based structure consisting of $q^3$ cubes, whose side length is $\delta_{cube} \equiv \delta_{subdom}$. Then, the sets ${\cal X}_n$, ${\cal C}_d$ and ${\cal E}_s$ are partitioned by the cube structure into $q^3$ subsets ${\cal X}_{n_k}$, ${\cal C}_{d_k}$ and ${\cal E}_{p_k}$, $k=1,2,\ldots,q^3$, where $n_k$, $d_k$ and $p_k$ are the number of points in the $k$-th cube.

\vskip 3pt 
This stage can be summarized in Algorithm \ref{alg1}.  
\vskip 3pt 

\begin{algorithm}[ht!]                   
\caption{Cube-partition structure}          
\label{alg1}                          
\begin{algorithmic}[1]                   
\FOR{each cube $k = [u,v,w]$, $u,v,w=1,2,\ldots,q$}
\STATE partition and count the number of points
\STATE \hspace{1.cm} $n_k = n_{u,v,w}$ (nodes)
\STATE \hspace{1.cm} $d_k = d_{u,v,w}$ (subdomain points)
\STATE \hspace{1.cm} $p_k = p_{u,v,w}$ (evaluation points);
\RETURN $(n_k; {\cal X}_{n_k})  \wedge (d_k;{\cal C}_{d_k}) \wedge (p_k;{\cal E}_{p_k})$
\ENDFOR
\end{algorithmic}
\end{algorithm}

\vskip3pt
\noindent {\tt Stage 5.} In order to identify the cubes to be examined in the searching procedure, we adopt the following rule which is composed of three steps: 
\begin{enumerate}
\item[(1)] the cube side $\delta_{cube}$ is chosen equal to the subdomain radius $\delta_{subdom}$, i.e. $\delta_{cube} \equiv \delta_{subdom}$, and the ratio between these quantities is denoted by $i^* = \delta_{subdom}/\delta_{cube}$;
\item[(2)] the value $i^*$ provides the number $j^*$ of cubes to be examined for each point by the rule $j^* = (2i^*+1)^3$, which obviously here gives $ j^*=27$. In practice, this means that the search of the nearby points is limited at most to twenty-seven cubes: the cube on which the considered point lies, and the twenty-six neighboring cubes; 
\item[(3)] for each cube $k=[u,v,w]$, $u,v,w=1,2,\ldots,q$, a cube-partition searching procedure is considered, examining the points from the cube $[u-i^*,v-i^*,w-i^*]$ to the cube $[u+i^*,v+i^*,w+i^*]$. For the points of the first and last cubes (those close to the boundary of the unit cube), we reduce the total number of cubes to be examined, setting $u-i^* = 1$ and/or $v-i^* = 1$ and/or $w-i^* = 1$ (when $u-i^* < 1$ and/or $v-i^* < 1$ and/or $w-i^* < 1$) and $u+i^* = q$ and/or $v+i^* = q$ and/or $w+i^* = q$ (when $u+i^* > q$ and/or $v+i^* > q$ and/or $w+i^* > q$). 
\end{enumerate}

Then, after defining which and how many cubes are to be examined, the cube-partition searching procedure (see Algorithm \ref{alg2}) is applied:
\begin{itemize}
	\item for each subdomain point of ${\cal C}_{d_k}$, $k=1,2,\ldots,q^3$, to determine all nodes belonging to a subdomain. The number of nodes of the subdomain centred at $(\bar{x}_j,\bar{y}_j,\bar{z}_j)$ is counted and stored in $\bar{n}_j$, $j=1,2,\ldots,d$; 
	\item for each evaluation point of ${\cal E}_{p_k}$, $k=1,2,\ldots,q^3$, in order to find all those belonging to a subdomain of centre $(\bar{x}_i,\bar{y}_i,\bar{z}_i)$ and radius $\delta_{subdom}$. The number of subdomains containing the $i$-th evaluation point is counted and stored in $r_i$, $i=1,2,\ldots,s$. 
\end{itemize}


\begin{algorithm}[ht!]               
\caption{Cube-partition searching procedure}          
\label{alg2}                          
\begin{algorithmic}[1]                   
\FOR{$w = 1,2,\ldots,q$}
\FOR{$v = 1,2,\ldots,q$}
\FOR{$u = 1,2,\ldots,q$}
\STATE set $[first_x,first_y,first_z] = [u - i^*,v - i^*,w - i^*]$
\STATE \hspace{0.5cm} $[last_x,last_y,last_z] = [u + i^*,v + i^*,w + i^*]$
\IF{$first_x < 1$ and/or $first_y < 1$ and/or $first_z < 1$}
\STATE set $first_x = 1$ and/or $first_y = 1$ and/or $first_z = 1$
\ENDIF
\IF{$last_x > q$ and/or $last_y > q$ and/or $last_z > q$}
\STATE set $last_x = q$ and/or $last_y = q$ and/or $last_z = q$
\ENDIF
\FOR{$h = subdom\_bp_{u,v,w}, \ldots, subdom\_ep_{u,v,w}$}
\STATE set $\bar{n}_h = 0$
\FOR{$k = first_z, \ldots, last_z$}
\FOR{$j = first_y, \ldots, last_y$}
\FOR{$i = first_x, \ldots, last_x$}
\FOR{$r = bp_{i,j,k}, \ldots, ep_{i,j,k}$}
\IF{$(x_r,y_r,z_r) \in I_h((\bar{x},\bar{y},\bar{z});\delta_{subdom})$}
\STATE set $\bar{n}_h = \bar{n}_h + 1$
\STATE \hspace{0.5cm} $STORE_{h,\bar{n}_h}(x_r,y_r,z_r,f_r)$
\ENDIF
\ENDFOR
\ENDFOR
\ENDFOR
\ENDFOR
\RETURN $(x,y,z)\in I_h((\bar{x},\bar{y},\bar{z});\delta_{subdom})$
\ENDFOR

\FOR{$h = eval\_bp_{u,v,w}, \ldots, eval\_ep_{u,v,w}$}
\STATE set $r_h = 0$
\FOR{$k = first_z, \ldots, last_z$}
\FOR{$j = first_y, \ldots, last_y$}
\FOR{$i = first_x, \ldots, last_x$}
\FOR{$r = subdom\_bp_{i,j,k} \ldots, subdom\_ep_{i,j,k}$}
\IF{$(\tilde{x}_r,\tilde{y}_r,\tilde{z}_r) \in I_h((\bar{x},\bar{y},\bar{z});\delta_{subdom})$}
\STATE set $r_h = r_h + 1$
\STATE \hspace{0.5cm} $STORE_{h,r_h}(\tilde{x}_r,\tilde{y}_r,\tilde{z}_r)$
\ENDIF
\ENDFOR
\ENDFOR
\ENDFOR
\ENDFOR
\RETURN $(\tilde{x},\tilde{y},\tilde{z})\in I_h((\bar{x},\bar{y},\bar{z});\delta_{subdom})$
\ENDFOR

\ENDFOR
\ENDFOR
\ENDFOR
\end{algorithmic}
\end{algorithm}

\vskip 3pt 

\noindent {\tt Stage 6.} A local interpolant $R_j$, $j=1,2, \ldots, d$, is found for each subdomain point.
\vskip 3pt 
\noindent {\tt Stage 7.} A local approximant $R_j(x,y,z)$ and a weight function $W_j(x,y,z)$, $j=1,2,\ldots,d$, is found for each evaluation point.
\vskip 3pt 
\noindent {\tt Stage 8.} Applying the global interpolant (\ref{pui}), one can find approximated values computed at any evaluation point $(\tilde{x},\tilde{y},\tilde{z}) \in {\cal E}_s$.


\subsection{Complexity analysis} \label{comp_anal}

The algorithm is based on the construction of a cube-partition searching procedure. It enables us to efficiently determine all points belonging to each subdomain $\Omega_j$, $j=1,2,\ldots,d$, so that we can compute local RBF interpolants to be used in the partition of unity scheme. Assuming that the covering $\{\Omega_j\}_{j=1}^{d}$ is regular and local and the set ${\cal X}_n$ of data points is quasi-uniform, we analyze the complexity of this code. 

The cube-partition algorithm involves the use of the standard \textsl{quicksort} routine, which requires on average a time complexity ${\cal O}(M\log M)$, where $M$ is the number of points to be sorted. Specifically, we have a distribution phase consisting of building the data structure, in which the computational cost has order: ${\cal O}(n\log n)$ for the sorting of all $n$ nodes and ${\cal O}(s\log s)$ for the sorting of all $s$ evaluation points in \texttt{Stage 1}. Then, \texttt{in Stage 3} the \textsl{quicksort} routine is repeatedly used with respect to different directions considering a reduced number of points (see Subsections \ref{cube-proc}--\ref{cube-algor}). Since the number of centres in each subdomain $\Omega_j$ is bounded by a constant (see Definition \ref{defpr}), we need ${\cal O}(1)$ space and time for each subdomain to solve the local RBF interpolation problems. In fact, in order to obtain the local RBF interpolants, we have to solve $d$ linear systems of (relatively) small sizes, i.e. $\bar{n}_j \times \bar{n}_j$, with $\bar{n}_j << n$,  thus requiring a constant running time ${\cal O}(\bar{n}_j^3)$, $j=1,2,\ldots,d$, for each subdomain (see \texttt{Stage 6}). Then, in \texttt{Stage 5, 7} and \texttt{8} we also need a cost of $r_k \cdot {\cal O}(\bar{n}_j)$, $j=1,2,\ldots,d$, $k=1,2,\ldots,s$, for the $k$-th evaluation point of ${\cal E}_s$; in other words, we have a constant time to get the value of the global fit \eqref{pui}. Finally, the algorithm requires $4n$, $4d$ and $4s$ storage requirements for the data, and $\bar{n}_j$, $j=1,2,\ldots, d$, locations for the coefficients of each local RBF interpolant.


\section{Numerical experiments} \label{num_res}
In this section we present a few numerical tests to show performance of the cube-partition algorithm, numerically analyzing efficiency and accuracy of the local interpolation scheme on some sets of scattered data. The code is implemented in C/C++ language, while numerical results are carried out on a Intel Core i7-4500U 1.8 GHz processor. In the experiments we consider a node distribution containing $n=(2^k + 1)^3$, $k=4,5,6$, uniformly random Halton nodes generated by using the MATLAB program \texttt{haltonseq.m} (see \cite{Fasshauer07}). The cube-partition algorithm is run considering $d=8^{k-1}$, $k=4,5,6$, subdomain points and $s = 11^3=1331$ evaluation (or grid) points, which are contained in the unit cube $\Omega = [0,1]^3$. Here, for the global interpolant (\ref{pui}) we use Shepard's weight (\ref{sh_w}).


The performance of the interpolation algorithm is verified taking the data values by the following two trivariate Franke's test functions (see, e.g., \cite{Lazzaro02,Renka88a})
\begin{eqnarray}
\hskip0.7cm f_1(x,y,z)&=&\frac{3}{4}\exp\left[-\frac{(9x-2)^2+(9y-2)^2+(9z-2)^2}{4}\right]  \nonumber \\
&+&\frac{3}{4} \exp\left[-\frac{(9x+1)^2}{49}-\frac{9y+1}{10}-\frac{9z+1}{10}\right] \nonumber \\
&+&\frac{1}{2} \exp\left[-\frac{(9x-7)^2+(9y-3)^2+(9z-5)^2}{4}\right] \nonumber \\
&-&\frac{1}{5} \exp\left[-(9x-4)^2-(9y-7)^2-(9z-5)^2\right],\nonumber \\ \nonumber \\
f_2(x,y,z)&=&\frac{\left(1.25+\cos(5.4y)\right)\cos(6z)}{6+6\left(3x-1\right)^2}, \nonumber
\end{eqnarray}
and using Gaussian $C^{\infty}$ (G), Mat$\acute{\text{e}}$rn $C^4$ (M4) and  Wendland $C^4$ (W4) as local RBF interpolants
\begin{equation}
\left.
\begin{array}{rcllc}
\phi(r) & = & {\rm e}^{-\alpha^2 r^2}, &                      & \hspace{1.cm} {\rm G} \nonumber \\
\phi(r) & = & \displaystyle{{\rm e}^{-\epsilon r} (\epsilon^2r^2+3\epsilon r+3)},&   & \hspace{1.cm} \mbox{{\rm M4}} \nonumber\\
\phi(r) & = & \displaystyle{\left(1-\delta r\right)_+^6(35\delta^2r^2+18\delta r+3)},&   & \hspace{1.cm} \mbox{{\rm W4}} \nonumber
\end{array}
\right.
\end{equation}
where $\alpha, \epsilon, \delta \in \RR^+$ are the \textsl{shape parameters}, $r=||\cdot||_2$ is the Euclidean distance, and $(\cdot)_+$ denotes the truncated power function. Note that Gaussian $C^{\infty}$ and Mat$\acute{\text{e}}$rn $C^4$ are globally supported basis functions, whereas Wendland $C^4$ is a compactly supported one (see \cite{Wendland05}).

Some information about the execution of the interpolation algorithm described in Section \ref{PUM_ALG} are reported in Table \ref{time}, namely the number $q^3$ of partitions in cubes of the domain and the CPU times (in seconds) obtained by running the cube-partition algorithm. Moreover, since we are interested in pointing out the effectiveness of the proposed algorithm, in Table \ref{time} we also show CPU times obtained by using the same interpolation method, but without partitioning the domain $\Omega$ in cubes and, accordingly, without considering the corresponding searching procedure. This analysis emphasizes that the use of a cube structure gives a considerable saving of time, mainly when the number of points to be handled becomes quite a lot large.

\begin{table}[ht!]
		\begin{center}
			\begin{tabular}{|cc|cc|c|} 
			 \hline
			 \rule[0mm]{0mm}{3ex}
			  	$n$    & $d$ & $q^3$ & \textbf{$t_{cube}$} & $t_{no-cube}$\\
				\hline
				\rule[0mm]{0mm}{3ex}
				$4913$     & $512$    & $6^3$  & $\textbf{1.1}$   & $1.4$     \\ 
				\rule[0mm]{0mm}{3ex}
				$35937$    & $4096$   & $12^3$ & $\textbf{7.9}$  & $15.5$    \\ 
			  \rule[0mm]{0mm}{3ex}
			  $274625$   & $32768$  & $23^3$ & $\textbf{62.7}$ & $525.0$  \\ 
			  \hline 
			\end{tabular}
		\end{center}
			\caption{Number of partitions in cubes and CPU times (in seconds) obtained by running the cube-partition algorithm ($t_{cube}$), and the corresponding one without a cube structure ($t_{no-cube}$).}
			\label{time}
	\end{table} 

Analyzing the performance of the algorithm, we observe that the cube-partition searching procedure turns out to be powerful and efficient, because CPU times reported in Table \ref{time} are mainly due to solution of $d$ linear systems having matrices with a relatively large number of entries, usually more than a hundred. 

\begin{figure}[ht!]
\begin{center}
\begin{minipage}{70.mm}
\includegraphics[width=7.cm]{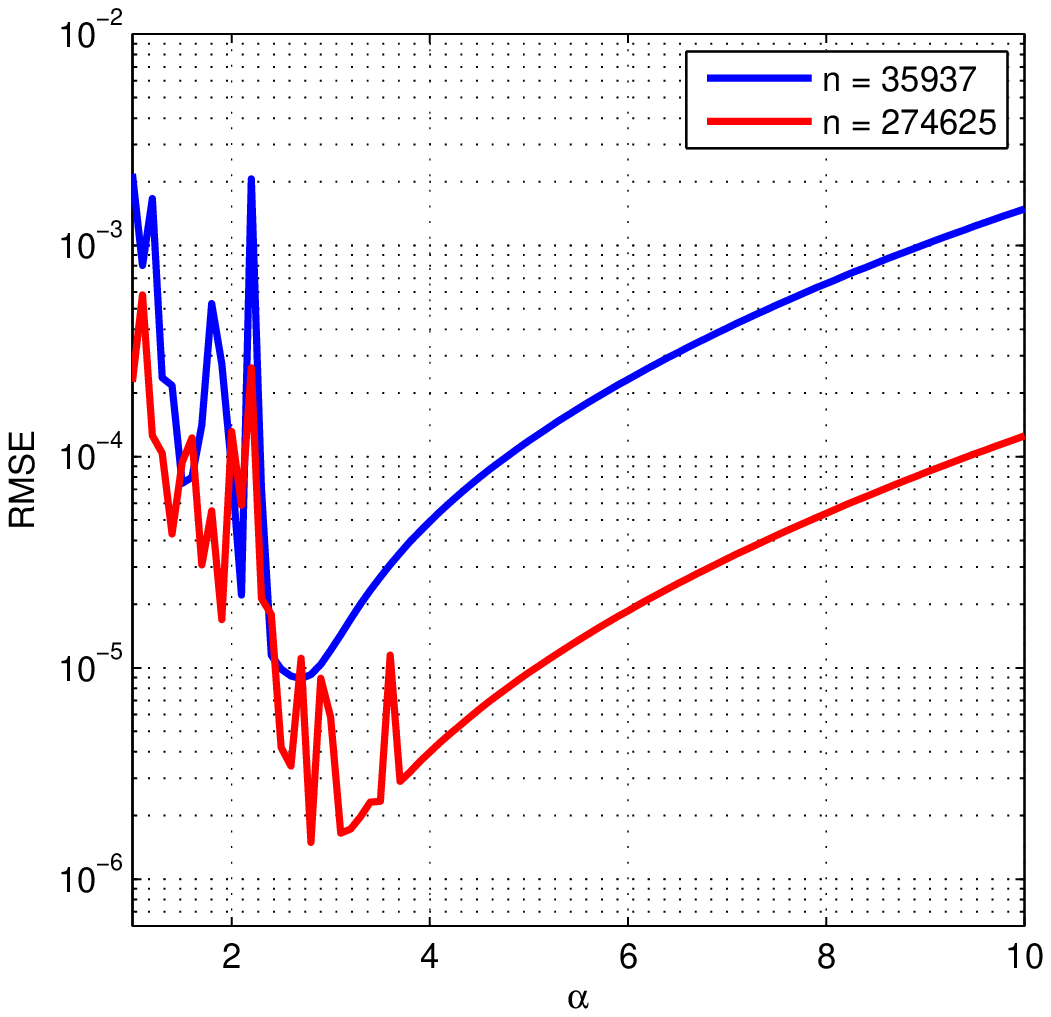}
\centerline{G -- $f_1$}
\end{minipage}\hskip -1.2cm
\begin{minipage}{70.mm}
\includegraphics[width=7.cm]{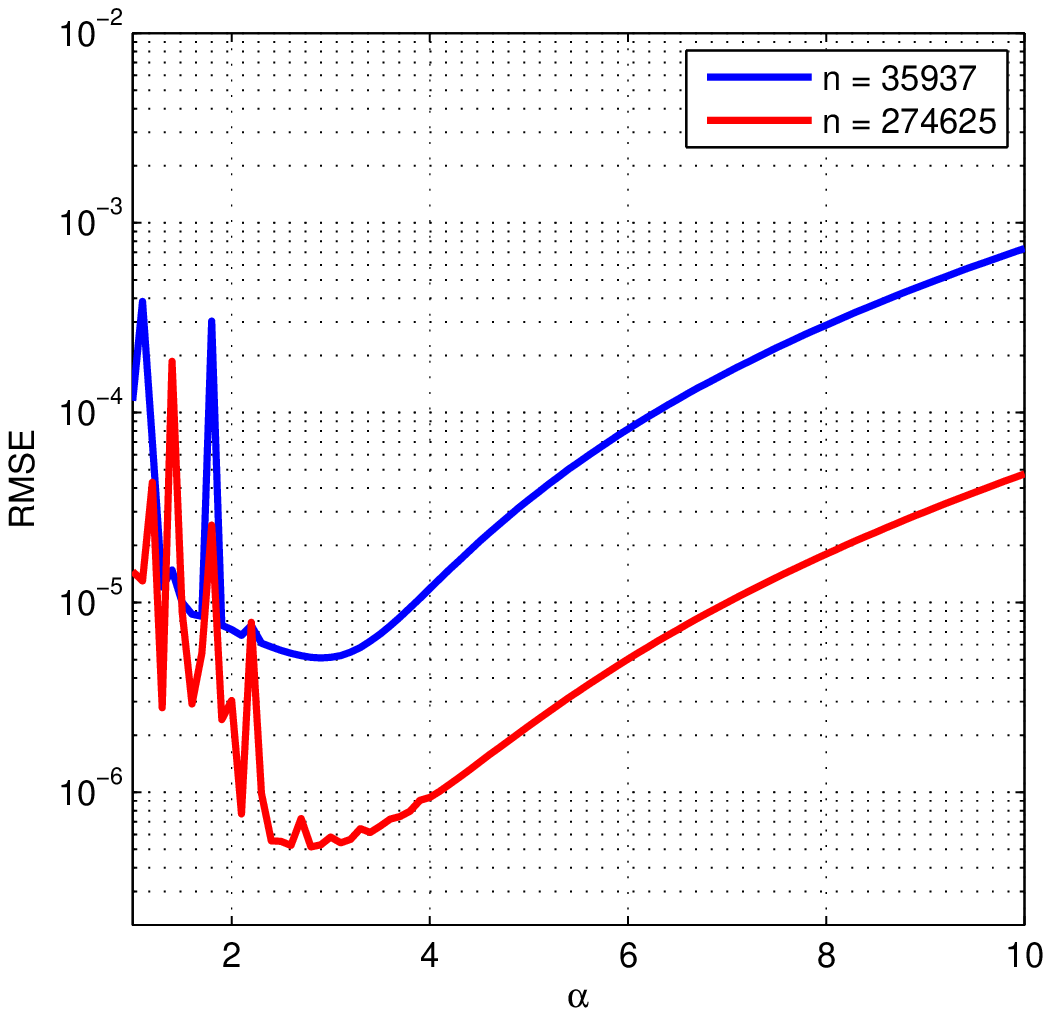}
\centerline{G -- $f_2$}
\end{minipage}\\
\begin{minipage}{70.mm}
\includegraphics[width=7.cm]{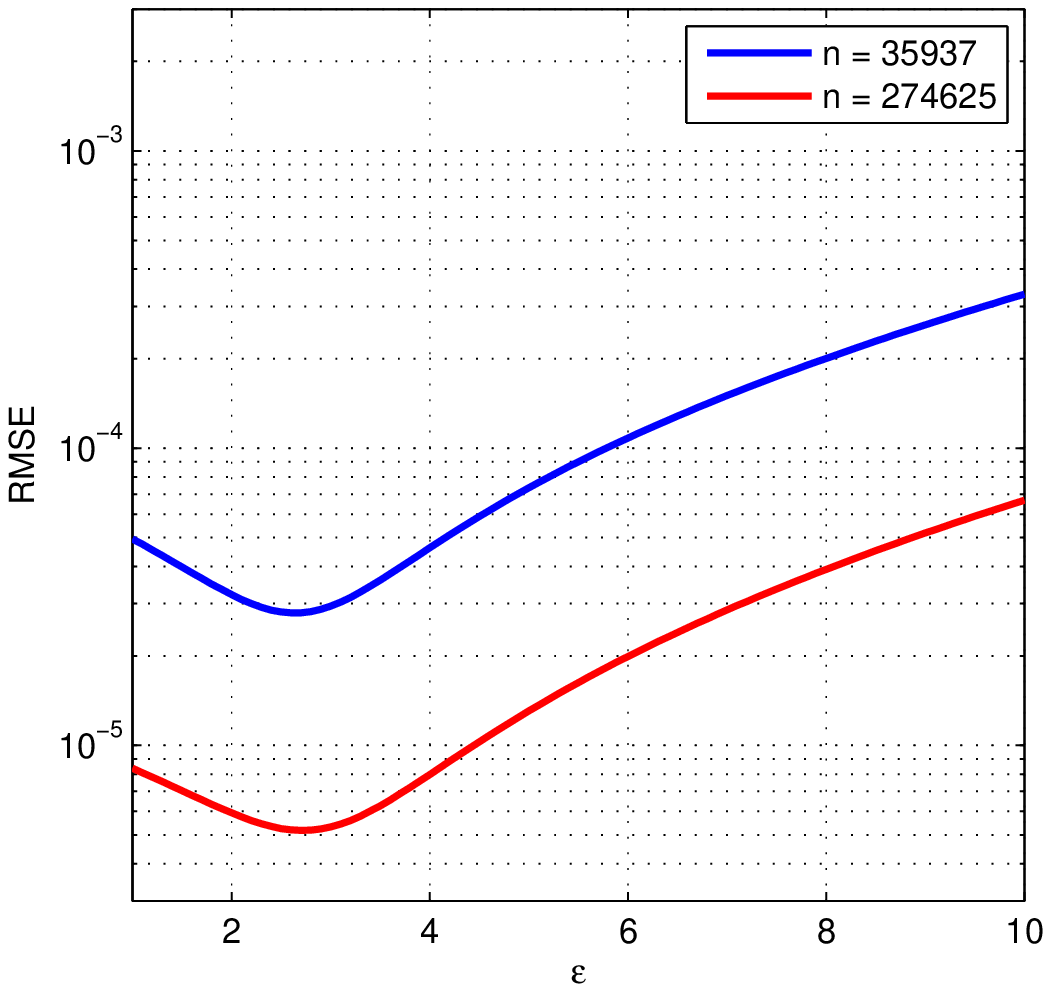}
\centerline{M4 -- $f_1$}
\end{minipage}\hskip -1.2cm
\begin{minipage}{70.mm}
\includegraphics[width=7.cm]{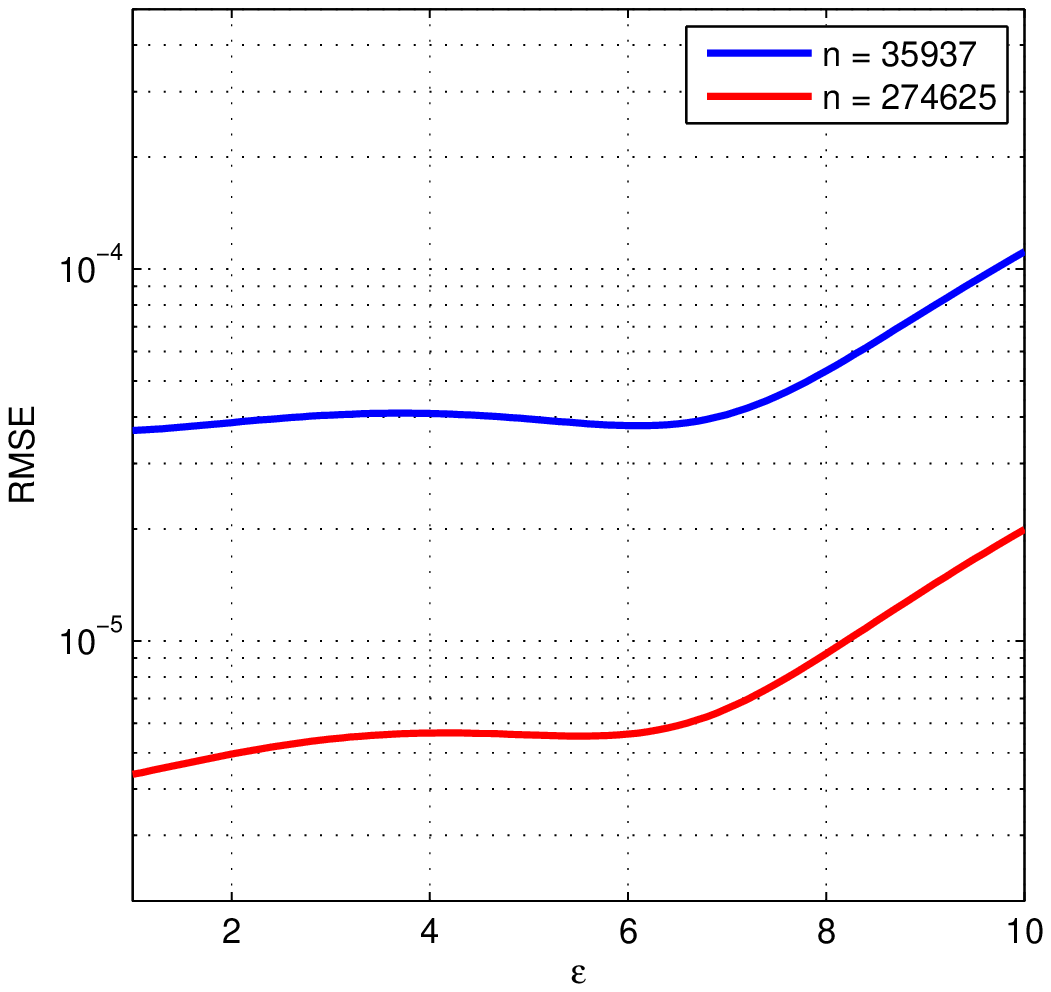}
\centerline{M4 -- $f_2$}
\end{minipage}\\
\begin{minipage}{70.mm}
\includegraphics[width=7.cm]{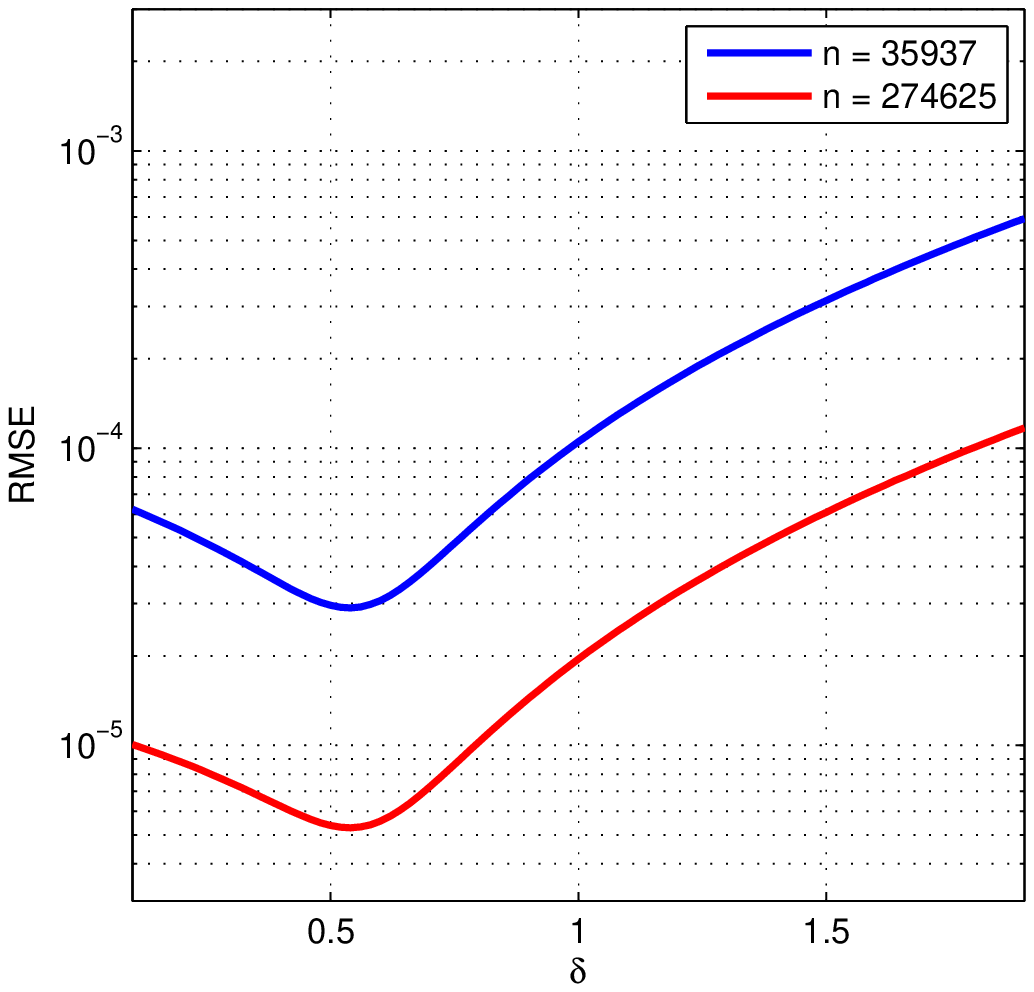}
\centerline{W4 -- $f_1$}
\end{minipage}\hskip -1.2cm
\begin{minipage}{70.mm}
\includegraphics[width=7.cm]{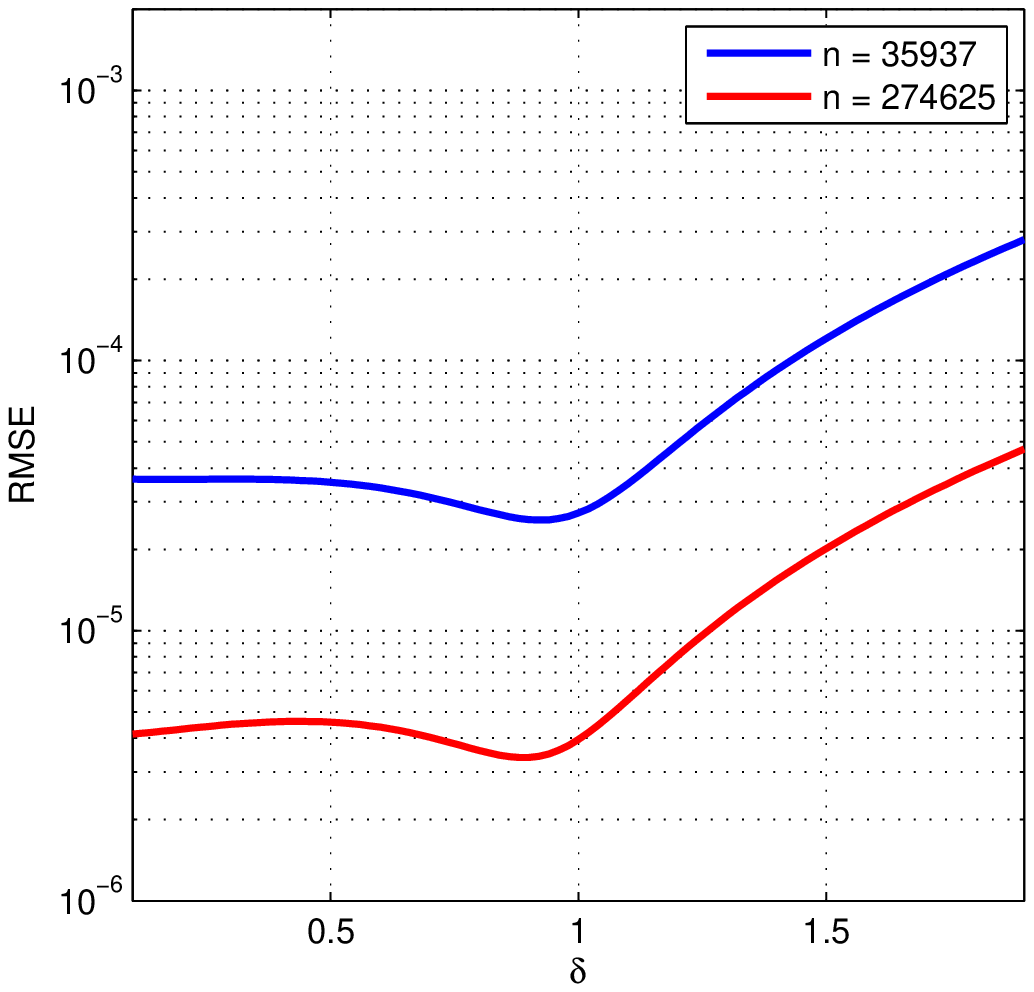}
\centerline{W4 -- $f_2$}
\end{minipage}
\end{center}
\caption{RMSEs obtained by varying the shape parameters.}
\label{shape_f1}
\end{figure}

Now, in order to investigate accuracy of the method, we compute the root mean square error (RMSE), whose formula is
\begin{equation}
	RMSE = \sqrt{\frac{1}{s}\sum_{i=1}^{s} |f(\boldsymbol{x}_i) - {\cal I}(\boldsymbol{x}_i)|^2}, \nonumber
\end{equation}
analyzing its behavior by varying the values of the shape parameters for Gaussian, Mat$\acute{\text{e}}$rn and Wendland functions (see Figure \ref{shape_f1}). These graphs allow us to find the optimal values of $\alpha$, $\epsilon$ and $\delta$, i.e. those values for which we obtain the smallest RMSEs (see Tables \ref{tab_errors}--\ref{tab_errors_bis}). Note that each evaluation is carried out by choosing equispaced values of the shape parameters, taking $\alpha, \epsilon \in [1, 10]$ and $\delta \in [0.1, 1.9]$. Analyzing error tables and graphs, we can see that Mat$\acute{\text{e}}$rn and Wendland functions have a greater stability than RBF Gaussian, but the latter gives us a greater accuracy although its interpolation matrices might be subject to ill-conditioning problems for small values of $\alpha$. This behavior is what we expect from theoretical standpoint, but here it is validated by numerical tests. Moreover, we remark that several numerical experiments (not reported here for shortness) have been carried out using other test functions and the results show a uniform behavior.


\begin{table}[ht!]
		\begin{center}
			\begin{tabular}{|c|c|c|c|c|c|c|} \hline
			  & \multicolumn{2}{c|}{  \rule[-2mm]{0mm}{7mm}  G} & \multicolumn{2}{c|}{  \rule[-2mm]{0mm}{7mm}  M4} & \multicolumn{2}{c|}{  \rule[-2mm]{0mm}{7mm}  W4}  \\
			  \cline{2-7} \rule[-2mm]{0mm}{7mm}
			  $n$	& RMSE & $\alpha_{opt}$ & RMSE & $\epsilon_{opt}$ & RMSE & $\delta_{opt}$ \\
				\hline 
				\rule[0mm]{0mm}{3ex}
				$\hskip-1pt 35937$ & $8.8797{\rm E}-6$  & $2.7$  & $2.7905{\rm E}-5$  & $2.6$  & $2.9041{\rm E}-5$  & $0.54$ \\
				\rule[0mm]{0mm}{3ex}
				$274625$           & $1.4928{\rm E}-6$  & $2.8$  & $5.1734{\rm E}-6$  & $2.7$  & $5.2847{\rm E}-6$  & $0.54$ \\
			  \hline 
			\end{tabular}
		\end{center}
			\caption{RMSEs obtained by using optimal values of $\alpha$, $\epsilon$ and $\delta$ for $f_1$.}
			\label{tab_errors}
	\end{table}

\begin{table}[ht!]
		\begin{center}
			\begin{tabular}{|c|c|c|c|c|c|c|} \hline
			  & \multicolumn{2}{c|}{  \rule[-2mm]{0mm}{7mm}  G} & \multicolumn{2}{c|}{  \rule[-2mm]{0mm}{7mm}  M4} & \multicolumn{2}{c|}{  \rule[-2mm]{0mm}{7mm}  W4}  \\
			  \cline{2-7} \rule[-2mm]{0mm}{7mm}
			  $n$	& RMSE & $\alpha_{opt}$ & RMSE & $\epsilon_{opt}$ & RMSE & $\delta_{opt}$ \\
				\hline 
				\rule[0mm]{0mm}{3ex}
				$\hskip-1pt 35937$ & $5.1013{\rm E}-6$  & $2.9$  & $3.6761{\rm E}-5$  & $1.0$  & $2.5677{\rm E}-5$  & $0.92$ \\
				\rule[0mm]{0mm}{3ex}
				$274625$           & $5.1446{\rm E}-7$  & $2.8$  & $4.3760{\rm E}-6$  & $1.0$  & $3.3941{\rm E}-6$  & $0.88$ \\
			  \hline 
			\end{tabular}
		\end{center}
			\caption{RMSEs obtained by using optimal values of $\alpha$, $\epsilon$ and $\delta$ for $f_2$.}
			\label{tab_errors_bis}
	\end{table}

Finally, to show that the CPU times in Table \ref{time} essentially depend on the size of interpolation matrices, we repeat numerical tests fixing a maximum number (i.e., $m_i=m_{max}$, $i=1,2,\ldots,d$) of nodes for each subdomain, namely only considering the $m_{max}$ nodes closest to the subdomain centres. In fact, for example, taking $m_{max}=50,70$ (and also $m_{max}$ not fixed) and denoting by $t_{cube}^{m_{max}}$ the corresponding execution times, we get a significant reduction of times, since $t_{cube}^{50} = 0.5$ and $t_{cube}^{70} = 0.6$ for $n = 4913$, $t_{cube}^{50} = 1.9$ and $t_{cube}^{70} = 3.4$ for $n = 35937$, while $t_{cube}^{50} = 14.2$ and $t_{cube}^{70} = 28.1$ for $n=274625$ (see Table \ref{time} for a comparison). Nevertheless, this reduction expressed in terms of CPU times is paid, in general, only with a slight loss of accuracy, since the behavior of RMSEs is similar to that shown in Figure \ref{shape_f1}. 


In conclusion, in Table \ref{tab_comp_1} we also report the RMSEs obtained by applying the cube-partition algorithm on sets of grid points.

\begin{table}[ht!]
		\begin{center}
			\begin{tabular}{|c|c|c|c|c|c|c|} \hline
			  $n$ & \multicolumn{2}{c|}{  \rule[-2mm]{0mm}{7mm}  $35937$} & \multicolumn{2}{c|}{  \rule[-2mm]{0mm}{7mm}  $274625$}  \\
			  \cline{1-5} \rule[-2mm]{0mm}{7mm}
			        & $f_1$ & $f_2$ & $f_1$ & $f_2$  \\
				\hline 
				\rule[0mm]{0mm}{3ex} 
				      G     & $2.4327{\rm E}-6$  & $1.5521{\rm E}-7$      & $2.6580{\rm E}-7$  & $1.3038{\rm E}-8$   \\
	  $\alpha_{opt}$  & $3.4$                & $3.1$                    & $4.4$                & $2.7$   \\
				\rule[0mm]{0mm}{3ex}
				     M4     & $1.3052{\rm E}-5$  & $3.0937{\rm E}-6$      & $2.9642{\rm E}-6$  & $6.2521{\rm E}-7$   \\
	  $\epsilon_{opt}$  & $2.0$              & $1.5$                    & $4.3$                & $2.2$   \\
				\rule[0mm]{0mm}{3ex}
				     W4     & $1.2938{\rm E}-5$  & $2.8769{\rm E}-6$      & $2.6711{\rm E}-6$  & $5.9111{\rm E}-7$   \\
	  $\delta_{opt}$  & $0.48$               & $0.50$                    & $0.86$                & $0.46$   \\
			  \hline 
			\end{tabular}
		\end{center}
			\caption{RMSEs computed on grid points.}
			\label{tab_comp_1}
	\end{table}


\section{Conclusions and future work} \label{concl}
In this paper we propose a new local interpolation algorithm for trivariate interpolation of scattered data points. It is based on the construction of a partition of the domain in cubes, enabling us to optimally implement a cube-partition searching procedure in order to efficiently detect the nodes belonging to each subdomain of the partition of unity method. This technique works well and quickly also when the amount of data to be interpolated is very large. Moreover, the proposed algorithm is flexible, since different choices of local interpolants are allowable, and completely automatic. 

As regards research and future work we are interested in refining the cube algorithm adopting suitable data structures like kd-trees and range trees, connecting these data structures with the special partition of the domain in cubes. Moreover, we are going to extend the proposed algorithm to higher dimensions. Then, even though the choice of low-order basis functions such as Mat$\acute{\text{e}}$rn and Wendland functions gives a good trade-off between stability and accuracy, we are still considering the need of dealing with the ill-conditioning problem of high-order basis functions. On the one hand, we might consider suitable preconditioning techniques for RBF interpolation matrices as already done in \cite{Cavoretto12b} for RBF collocation matrices; on the other hand, one could study alternative stategies to have a stable evaluation of interpolants via Hilbert-Schmidt SVD as in \cite{Cavoretto14d,Fasshauer12}, or new stable bases as in \cite{DeMarchi13,Pazouki11}. 






\end{document}